\documentclass[a4paper,12pt,pagesize,headings=small]{scrartcl}
\pdfoutput=1
\usepackage[utf8]{inputenc}
\usepackage[english,ngerman]{babel}
\usepackage[T1]{fontenc}
\usepackage{lmodern}
\usepackage{microtype}
\usepackage{amsmath}
\usepackage{amssymb}
\usepackage{hyperref}
\usepackage[all]{xy}
\usepackage[hyperref,amsmath,thmmarks]{ntheorem}
\usepackage{tikz}
\usetikzlibrary{matrix,arrows}

\newcommand{\Z}{\mathbb Z}

\newcommand{\F}{\mathbb F}

\newcommand{\cA}{\mathcal A}

\newcommand{\cW}{\mathcal W}
\newcommand{\cQ}{\mathcal Q}

\newcommand{\cU}{\mathcal U}

\DeclareMathOperator{\Gl}{GL}
\DeclareMathOperator{\E}{E}
\DeclareMathOperator{\Sl}{SL}
\DeclareMathOperator{\St}{St}

\DeclareMathOperator{\Flag}{Flag}

\DeclareMathOperator{\rank}{rank}

\DeclareMathOperator{\lk}{lk}

\newcommand{\coloneq}{\mathrel{\mathop :}=}
\newcommand{\alter}[2]{\left\{\genfrac{}{}{0pt}{}{#1}{#2} \right\}}

\theoremstyle{plain}
\newtheorem{Definition}{Definition}
\newtheorem{Theorem}[Definition]{Theorem}
\newtheorem{Lemma}[Definition]{Lemma}
\newtheorem{Proposition}[Definition]{Proposition}
\newtheorem{Corollary}[Definition]{Corollary}

\theoremclass{Theorem}
\theoremstyle{break}

\theoremstyle{nonumberplain}
\theorembodyfont{\rm}
\newtheorem{Remark}{Remark}
\theoremsymbol{\ensuremath\Box}
\newtheorem{Proof}{Proof}

\addtokomafont{paragraph}{\rm\bf}
\addtokomafont{subparagraph}{\rm\it}

\title{On Wagoner complexes}
\author{Jan Essert}
\date{\today}

\xdefinecolor{myblue}{rgb}{0.1,0.1,0.6}
\xdefinecolor{mygreen}{rgb}{0.1,0.4,0.3}

\hypersetup{
pdftitle={On Wagoner complexes},
pdfauthor={Jan Essert},
pdfkeywords={Wagoner complexes, Buildings, Group homology, Group theory, K-theory},
colorlinks=true,
linkcolor=myblue,
citecolor=mygreen,
urlcolor=myblue
}

\begin{document}
\selectlanguage{english}
\maketitle

\begin{abstract}
	Wagoner complexes are simplicial complexes associated to groups of Kac-Moody type. They admit interesting homotopy groups which are related to integral group homology if the root datum is of $2$-spherical type. We give a general definition of Wagoner complexes, exhibit some simple properties and calculate low dimensional homotopy groups. In addition, we give a definition of affine Wagoner complexes related to groups admitting a root datum with valuation.
\end{abstract}

\addtocounter{section}{-1}
\addtocounter{Definition}{-1}

\section{Introduction}

Wagoner complexes are certain simplicial complexes first constructed by Wagoner in \cite{Wag:BSK:73}. In his original construction, they were associated to the group $\Gl_n(R)$ over an arbitrary ring $R$. Wagoner used the homotopy groups of these complexes to give a tentative definition of higher algebraic $K$-theory. Later Anderson, Karoubi and Wagoner proved that this definition is equivalent to Quillen's definition of $K$-theory.

The complexes Wagoner defined are closely related to the building associated to the general linear group. Unlike the building, the Wagoner complex is not homotopy equivalent to a bouquet of spheres, instead it has interesting homotopy and homology groups. Wagoner complexes have already been generalised to reductive $F$-groups by Petzold in \cite{Pet:HDW:04}. There, rational homology groups of Wagoner complexes are used to construct rational representations of reductive groups.

In this paper, we define Wagoner complexes in the context of groups of Kac-Moody type (or groups with root data), which we believe to be the appropriate level of generality. Important examples of such groups are semi-simple Lie groups and algebraic groups as well as Kac-Moody groups. All definitions are very natural and the connection to the theory of buildings is clearly visible.

We exhibit interesting substructures in Wagoner complexes analogous to apartments in buildings and present some general properties of Wagoner complexes. The main result is the following characterisation of low dimensional homotopy groups of Wagoner complexes in the case where the associated root datum is of $2$-spherical type.

\begin{Theorem}
    Let $G$ be a group of Kac-Moody type and assume that its root datum is of $2$-spherical type. Let $G^\dagger$ be the little projective group and let $\hat G$ be the associated Steinberg group in the sense of Tits. Let $\cW(G)$ be the associated Wagoner complex.
    
    Then, in the generic case, we have
        \begin{align*}
                \pi_0(|\cW(G)|) &\cong G/G^\dagger. \\
                \intertext{If $G$ is of non-spherical type or $\rank(G)\geq 3$, we have}
                \pi_1(|\cW(G)|) &\cong \ker(\hat G\twoheadrightarrow G^\dagger).
        \end{align*}
        By a recent result by Caprace, this means $\pi_1(|\cW(G)|)\cong H_2(G^\dagger;\Z)$ except over very small fields.
\end{Theorem}
For the precise statement, see Theorem \ref{th:main_result} and Corollary \ref{cor:main_result}.

In addition to this result, we also make the connection to algebraic $K$-theory explicit. Finally, we give a definition of affine Wagoner complexes associated to groups admitting root data with valuations.

The structure of this text is as follows. In section \ref{sec:definition}, we define Wagoner complexes. We construct apartments and a projection onto the building for spherical Wagoner complexes in section \ref{sec:simple_properties}. In section \ref{sec:roots}, we discuss simple properties of pairs of roots. The definition of the Steinberg group in the sense of Tits is given in section \ref{sec:steinberg}. We prove the main result in section \ref{sec:homotopy} and show its relation to algebraic $K$-theory in section \ref{sec:ktheory}. Finally, we give the construction of affine Wagoner complexes in section \ref{sec:affine}.

We would like to thank Linus Kramer for the introduction to the topic of Wagoner complexes and countless useful discussions. In addition, we are very grateful for suggestions and answers by Pierre-Emmanuel Caprace and Richard Weiss. We are indebted to Guntram Hainke for very good advice on the construction in the non-spherical case.

The author was supported by the \selectlanguage{ngerman}\emph{Graduiertenkolleg: \glqq Analytische Topologie und Me\-ta\-geo\-me\-trie\grqq}\selectlanguage{english} while working on this topic. This work is part of the author's doctoral thesis at the \selectlanguage{ngerman}Universität Münster\selectlanguage{english}.

\section{The definition of Wagoner complexes}\label{sec:definition}
Throughout this text, let $(W,S)$ be a Coxeter system of rank at least 2 whose Coxeter diagram has no isolated nodes. Let $\Sigma$ be the associated Coxeter complex. Denote by $\Phi$ the set of roots of $\Sigma$ and let $G$ be a group with a twin root datum $(G,(U_\alpha)_{\alpha\in\Phi},H)$ of type $\Phi$. For the precise definition of groups with root data, also known as groups of Kac-Moody type, see the book \cite{Rem:GKM:02} by Rémy. Denote the subgroup of $G$ generated by all root groups, the \emph{little projective group}, by $G^\dagger$.

We write $R_s=\{ t\in \Sigma: s\subseteq t\}$ for the \emph{residue} of $s$ in $\Sigma$. We will say that a simplex $s$ is \emph{co-spherical} if the corresponding residue $R_s$ is of spherical type. The basic `building blocks' for Wagoner complexes will be the following groups.

\begin{Definition}
	For a co-spherical simplex $s\in\Sigma$ we define
	\[
		U_s \coloneq \langle U_\alpha : R_s \subseteq \alpha \in \Phi\rangle.
	\]
\end{Definition}

\noindent It is well known that there is a twin building of type $\Phi$ associated to the group $G$, see \cite[2.6.4]{Rem:GKM:02}. By the same reference, the positive half of the twin building is a thick Moufang building $\Delta$ on which $G$ acts strongly transitively. There is a standard apartment in the building which we identify with $\Sigma$, and there is a fundamental chamber $c_0\in\Sigma$. Using the Moufang property of the building, we obtain:

\begin{Remark}
    We have that $s\leq s'$ if and only if $R_s\supseteq R_{s'}$ which is equivalent to $U_s\leq U_{s'}$.
\end{Remark}

\noindent For every $s\in\Sigma$, the parabolic subgroup $P_s$ is the stabiliser subgroup of $s$ in $G$. We have $U_s\leq P_s$. 

\begin{Remark}
	The building $\Delta$ is isomorphic to the set of cosets $\{ gP_s : s\subseteq c_0, g\in G\}$ of standard parabolics partially ordered by reverse inclusion, see for example \cite[6.29]{AB:B:08}. We will henceforth identify the building with this partially ordered set.
\end{Remark}

\noindent Generalising the construction by Wagoner in \cite{Wag:BSK:73}, we define:

\begin{Definition}\label{def:wagonercomplex}
    The \emph{Wagoner complex $\cW(G)$} is the flag complex over the set of cosets $\{ gU_s : s\in \Sigma, s \text{ co-spherical}, g\in G\}$ partially ordered by inclusion, which we write as follows:
	\[
    \cW(G)\coloneq \Flag(\{ gU_s : s\in\Sigma, s\text{ co-spherical}, g\in G\}).
	\]
    We call a Wagoner complex \emph{spherical} if the root datum is of spherical type.
\end{Definition}

\section{Simple Properties of spherical Wagoner complexes}\label{sec:simple_properties}

Although we make no direct use of this, let us examine some interesting properties of spherical Wagoner complexes for further reference. Most of these results were already observed in \cite{Pet:HDW:04} or \cite{Wag:BSK:73}.

We restrict ourselves to spherical Wagoner complexes in this section, because the corresponding results are more interesting. We will comment on the general situation at the end of this section. For the rest of this section, assume that $\Sigma$ is a spherical Coxeter complex.

\begin{Definition}
	We define the \emph{standard apartment of $\cW(G)$} to be
	\[
		\cA \coloneq \Flag(\{ U_s : s \in \Sigma \}).
	\]
	Left translates of the standard apartment will be called \emph{apartments}.
\end{Definition}

\noindent A spherical Wagoner complex admits a natural projection onto the barycentric subdivision of the building, which can be described as the flag complex over the building.

\paragraph{Projection} Consider the following map of partially ordered sets
\begin{align*}
    \tilde p_\Delta:\{ gU_s: s\in\Sigma, g\in G\} &\rightarrow \Delta \\
    gU_s &\mapsto gs.
\end{align*}
This map is well-defined since $U_s\leq P_s$. In addition, it is easily seen to be order-preserving. In particular, it induces a surjective simplicial map $p_\Delta:\cW(G)\rightarrow \Flag(\Delta)$ onto the barycentric subdivision of the building $\Delta$.

It is easy to see that the projection $p_\Delta$ maps the standard apartment $\cA$ bijectively onto the barycentric subdivision of the standard apartment $\Flag(\Sigma)$ which justifies the terminology. In addition, the map $p_\Delta$ induces a surjective map on singular homology: by the theorem of Solomon and Tits, the singular homology of the building $\Delta$ is concentrated in top-dimension and it is generated by the apartments which are all images of apartments in $\cW(G)$.

\paragraph{Right $N$-action} There is a right action of the group $N$ on the vertices of the Wagoner complex $\{ gU_s: g\in G, s\in\Sigma\}$ by right multiplication:
\[
	gU_sn = gnU_{n^{-1}s}.
\]
Since $U_{c}\cap N=\{1\}$ for any chamber $c$, this action is free. It induces a free right $N$-action on $\cW(G)$ and we have $\cA n = n\cA$. As already observed in \cite{Pet:HDW:04}, the projection $|\cW(G)|\rightarrow |\cW(G)/N|$ is a fibration with discrete fibre $N$.
The projection $p_\Delta$ factors through $\cW(G)/N$, since
\[
	\tilde p_\Delta( gnU_{n^{-1}s} ) = gn(n^{-1}s) = gs = \tilde p_\Delta(gU_s).
\]
If we consider a complex analogous to the Wagoner complex but replacing the groups $U_s$ by the parabolic subgroups $P_s$, we obtain another simplicial complex with a right $N$-action. The quotient complex is isomorphic to the barycentric subdivision of the building. The following diagram summarises all these constructions.
\[
\begin{xy}
	\xymatrix{
	\cW(G) = \Flag(\{ gU_s : s\in\Sigma, g\in G\})\ar[d]^{\mod N}\ar[r]\ar[dr]^{p_\Delta} & \Flag(\{ gP_s : s\in\Sigma, g\in G\})\ar[d]^{\mod N} \\
		\cW(G)/N \cong \Flag(\{ gU_s : s\subseteq c_0, g\in G\}) \ar[r]& \Flag(\Delta)\cong \Flag(\{ gP_s : s\subseteq c_0, g\in G\})
	}
\end{xy}
\]

\paragraph{The general case} In the general case for non-spherical root data, the concept of apartments and the above projection still make sense. The projection is then onto the flag complex over all co-spherical simplices of $\Delta$, which is the simplicial complex underlying the \emph{Davis realisation} of the building $\Delta$ as in \cite[18.2]{Dav:GTC:08}. The Davis realisation is contractible, so the induced maps on homotopy and homology groups are necessarily trivial.

\section{Root intervals}\label{sec:roots}

Remember that a pair of roots $(\alpha,\beta)$ is called \emph{prenilpotent}, if the intersections $\alpha\cap\beta$ and $-\alpha\cap-\beta$ both contain a chamber. For any prenilpotent pair of roots $(\alpha,\beta)$, we write
\begin{align*}
    [\alpha,\beta] &\coloneq \{\gamma\in\Phi: (\alpha\cap\beta)\subseteq \gamma, (-\alpha\cap-\beta)\subseteq -\gamma \}, \\
    ]\alpha,\beta] &\coloneq [\alpha,\beta]\setminus\{\alpha\},\\
    [\alpha,\beta[ &\coloneq [\alpha,\beta]\setminus\{\beta\}.
\end{align*}
These sets of roots are always finite. For any of the three sets of roots
\[
    \Psi \in \{ [\alpha,\beta], \;]\alpha,\beta], [\alpha,\beta[\;\},
\]
we write
\[
U_{\Psi}\coloneq \langle U_\gamma : \gamma \in \Psi \rangle.
\]

\begin{Remark}
    For spherical Coxeter systems, a pair of roots $(\alpha,\beta)$ is prenilpotent if and only if $\alpha\neq-\beta$. In this case, the root interval can be expressed simply as
    \[
        [\alpha,\beta] = \{\gamma\in\Phi : (\alpha\cap\beta)\subseteq \gamma\}.
    \]
\end{Remark}

\noindent For any root $\alpha\in\Phi$, we denote by $s_\alpha$ the reflection at its wall $\partial \alpha$. In non-spherical Coxeter complexes, there are two possibilities for prenilpotent pairs of roots.

\begin{Lemma}[Proposition 3.165 and Lemma 8.45 in \cite{AB:B:08}]\label{lem:nested_or_finite_order}
    Let $(\alpha,\beta)$ be a prenilpotent pair of roots in $\Phi$. Then there are two possibilities
\begin{itemize}
    \item The product $s_\alpha s_\beta$ has infinite order. In this case, the roots $\alpha$ and $\beta$ are \emph{nested}, that is, we have $\alpha\subseteq\beta$ or $\beta\subseteq\alpha$.
    \item The product $s_\alpha s_\beta$ has finite order. Then every maximal simplex $s$ of the intersection of the walls $\partial\alpha\cap\partial\beta$ has codimension $2$ and the link $\lk(s)$ is spherical. We write $\bar\alpha=\alpha\cap\lk(s)$, $\bar \beta=\beta\cap\lk(s)$ for the induced roots in the link. The map given by
    \begin{align*}
        [\alpha,\beta] &\rightarrow [\bar\alpha,\bar\beta]\\
        \gamma & \mapsto \gamma\cap\lk(s)
    \end{align*}
    is a bijection.
\end{itemize}
\end{Lemma}

\noindent This can be used to prove the following characterisation.

\begin{Lemma}\label{lem:root_intervals}
    For any prenilpotent, non-nested pair of roots $(\alpha,\beta)$, there are either roots $\alpha'$ and $\beta'$ such that
    \[
        ]\alpha,\beta] = [\alpha',\beta],\qquad [\alpha,\beta[\; = [\alpha,\beta']
    \]
    or $]\alpha,\beta] = \{\beta\}$ and $[\alpha,\beta[\;=\{\alpha\}$.
\end{Lemma}

\begin{Proof}
    If the Coxeter complex $\Sigma$ is of rank two and finite, it is an ordinary $2m$-gon for some $m\geq 2$. Its roots are paths of length $m$. We denote the vertices of $\Sigma$ by $(x_i)_{i\in \Z/2m}$. This leads to an enumeration of the roots $(\alpha_i)_{i\in \Z/2m}$ such that every root $\alpha_i$ contains the vertices $\{x_i,x_{i+1},\ldots,x_{i+m}\}$. Then it is easy to see that
    \[
        [\alpha_0,\alpha_i]=\{\alpha_0,\alpha_1,\ldots,\alpha_i\}
    \]
    for $1\leq i\leq m-1$. This implies in particular $[\alpha_0,\alpha_1]=\{\alpha_0,\alpha_1\}$ and $]\alpha_0,\alpha_i] = [\alpha_1,\alpha_i]$ and $[\alpha_0,\alpha_i[\; = [\alpha_0,\alpha_{i-1}]$ for $i>1$. Since the numbering of vertices was arbitrary, the Lemma follows for rank two finite Coxeter complexes.

    If the rank of the Coxeter complex $\Sigma$ is greater than two, and if $(\alpha,\beta)$ is a prenilpotent, non-nested pair of roots, by Lemma \ref{lem:nested_or_finite_order}, we can reduce to the rank 2 case as follows. Fix a maximal simplex $s$ in the intersection of the walls $\partial\alpha\cap\partial\beta$ and write $\bar\alpha=\alpha\cap\lk(s)$, $\bar \beta=\beta\cap\lk(s)$ as above. Then consider the bijection
    \begin{align*}
        [\alpha,\beta] &\rightarrow [\bar\alpha,\bar\beta]\\
        \gamma & \mapsto \gamma\cap\lk(s)
    \end{align*}
    Now the lemma follows by applying the finite rank two case to $\lk(s)$ and using the bijection to transport the result to the ambient building.
\end{Proof}

\noindent Using this characterisation, we prove another lemma which describes the `shadows' of root group intervals in the groups $U_s$.

\begin{Lemma}\label{lem:stabiliser_intersect_interval}
    Let $(\alpha,\beta)$ be a prenilpotent, non-nested pair of roots. Let $s\in\alpha\cap\beta$ be a co-spherical simplex. Then one of the following holds.
    \begin{itemize}
        \item There is a prenilpotent, non-nested pair of roots $(\gamma,\delta)$ such that $R_s\subset \gamma\cap\delta$ and such that $U_s\cap U_{[\alpha,\beta]}\subseteq U_{[\gamma,\delta]}$.
        \item There is a root $\gamma$ containing $R_s$ such that $U_s\cap U_{[\alpha,\beta]}\subseteq U_\gamma$.
    \end{itemize}
\end{Lemma}

\begin{Proof}
    The set of chambers in $R_s$ is finite, and we denote the number of chambers in $R_s\setminus (\alpha\cap\beta)$ by $n$. If $n=0$, then $(\alpha,\beta)$ is already a pair of roots we are looking for. Otherwise, we will replace $U_{[\alpha,\beta]}$ iteratively by smaller groups as follows, always decreasing $n$ strictly.

    If $n>0$, then $R_s\not\subset \alpha$ or $R_s\not\subset\beta$. We fix a chamber $c\in R_s\cap\alpha\cap\beta$, and a chamber $d\in -\alpha\cap -\beta$. Denote by $e$ the projection of $d$ onto $R_s$. Since the roots $-\alpha$ and $-\beta$ are convex and $d\in-\alpha\cap-\beta$, we have $e\in-\alpha\cap R_s$ or $e\in-\beta\cap R_s$. By \cite[2.2.6]{Rem:GKM:02}, the minimal gallery from $c$ to $d$ through $e$ crosses one of the walls $\partial\alpha$ or $\partial\beta$ first among the walls of roots in $[\alpha,\beta]$. Assume that the gallery crosses the wall $\partial\alpha$ first. Then there is a chamber $c'\in-\alpha\cap R_s$ such that $c'\in\gamma$ for all $\gamma\in\;]\alpha,\beta]$.

    By \cite[Proposition 8.33]{AB:B:08}, we have $U_{[\alpha,\beta]}=U_\alpha U_{]\alpha,\beta]}$. Then
    \[
    U_s \cap U_{[\alpha,\beta]} \subseteq U_s \cap U_{]\alpha,\beta]}
    \]
    since every non-trivial element of $U_\alpha$ cannot fix $c'$, while all elements of $U_{]\alpha,\beta]}$ and of $U_s$ fix $c'$. Now, by Lemma \ref{lem:root_intervals}, either there is a root $\alpha'$ such that $U_{]\alpha,\beta]}=U_{[\alpha',\beta]}$ or $]\alpha,\beta] = \{\beta\}$.

    In the first case, the number of chambers $n'$ in $R_s\setminus(\alpha'\cap\beta)$ is strictly smaller than the number of chambers in $R_s\setminus(\alpha\cap\beta)$, since the first one is a subset of the second and does not contain $c'$. If $n'=0$, we are finished, otherwise we repeat the procedure.

    In the second case, we have $U_{]\alpha,\beta]}=U_{\beta}$. Then either $R_s\in\beta$ and we are done, or $U_{[\alpha,\beta]}\cap U_s \subseteq U_s \cap U_{\beta}=\{1\}$, in which case we can choose an arbitrary root $\gamma$ containing $R_s$.
\end{Proof}

\section{The Steinberg group}\label{sec:steinberg}

From now on, we will often require the following technical condition. Denote by $\Pi$ the set of simple roots of the root system $\Phi$. For roots $\alpha,\beta\in\Pi$, we let
\[
X_\alpha \coloneq \langle U_\alpha \cup U_{-\alpha}\rangle,\qquad X_{\alpha,\beta} \coloneq \langle X_\alpha \cup X_\beta \rangle.
\]
Following Caprace in \cite{Cap:O2S:07}, we consider the following condition, which excludes the possibility that $X_{\alpha,\beta}$ modulo its centre is isomorphic to some small finite groups.\let\storedtheequation=\theequation\renewcommand{\theequation}{Co${}^*$}
\begin{equation}\label{condition}
    X_{\alpha,\beta}/Z(X_{\alpha,\beta}) \not\cong B_2(2), G_2(2), G_2(3),{}^2F_4(2)\quad\text{ for all pairs }\{\alpha,\beta\}\subseteq \Pi.
\end{equation}

\noindent We will require this condition for most of the results later on. Now, following Tits in \cite{Tit:Uni:87}, we introduce the Steinberg group. The definition uses the concept of the direct limit of a system of groups as in \cite[\S 1]{Ser:Tre:80}.\renewcommand{\theequation}{\storedtheequation}

\begin{Definition}
        The \emph{Steinberg group $\hat G$} associated to the group $G$ is the direct limit of the system of groups $U_\alpha$ and $U_{[\alpha,\beta]}$ with their natural inclusions for all roots $\alpha$ and all prenilpotent pairs of roots $(\alpha,\beta)$.
\end{Definition}

\noindent There is an improved characterisation of the Steinberg group due to Caprace, where we can restrict ourselves to non-nested pairs of roots.

\begin{Theorem}[Proposition 3.9 in \cite{Cap:O2S:07}]
    Assume that the group $G$ is of $2$-spherical type and satisfies condition \eqref{condition}. Then the Steinberg group $\hat G$ is isomorphic to the direct limit of the system of groups $U_\alpha$ and $U_{[\alpha,\beta]}$ with their natural inclusions for all roots $\alpha$ and all prenilpotent, non-nested pairs of roots $(\alpha,\beta)$.
\end{Theorem}

\noindent There is a similar characterisation for the groups $U_c$.

\begin{Theorem}[Theorem 3.6 in \cite{Cap:O2S:07}]\label{th:chamber_case}
    Assume that the group $G$ is of $2$-spherical type and satisfies condition \eqref{condition}. For any chamber $c\in\Sigma$, the group $U_c$ is the direct limit of the system of groups $U_\alpha$ and $U_{[\alpha,\beta]}$ for all roots $\alpha$ containing $c$ and all prenilpotent, non-nested pairs of roots $(\alpha,\beta)$ containing $c$, respectively.
\end{Theorem}

\noindent We need this statement also for the groups $U_s$.

\begin{Lemma}\label{lem:g_sinjection}
    Assume again that $G$ is of $2$-spherical type and satisfies condition \eqref{condition}. For any co-spherical simplex $s\in\Sigma$, the group $U_s$ is the direct limit of the system of groups $U_\alpha$ and $U_{[\alpha,\beta]}$ for all roots $\alpha$ where $R_s\subset\alpha$ and all prenilpotent, non-nested pairs of roots $(\alpha,\beta)$ with $R_s\subset(\alpha\cap\beta)$, respectively.
\end{Lemma}

\begin{Proof}
    If $s=c$ is a chamber, this is Theorem \ref{th:chamber_case}. Now let $s\in\Sigma$ be any co-spherical simplex and fix a chamber $c$ that contains $s$. Consider the following systems of groups
    \begin{align*}
        \cU_c &\coloneq \{ U_\alpha : c\in\alpha\} \cup \{ U_{[\alpha,\beta]} : c \in \alpha\cap\beta,\, (\alpha,\beta) \text{ non-nested}\} \\
        \cU_s &\coloneq \{ U_\alpha : R_s\subset\alpha\} \cup \{ U_{[\alpha,\beta]} : R_s \subset \alpha\cap\beta,\, (\alpha,\beta) \text{ non-nested}\}
    \end{align*}
    with the canonical inclusions. Since $\cU_s\subseteq \cU_c$, by the universal property of $\varinjlim\cU_s$, we obtain a homomorphism
    \[
        \varphi:\varinjlim\cU_s \rightarrow \varinjlim\cU_c
    \]
    which is the identity on all groups $U\in \cU_s$. We have the following diagram:
    \begin{center}
    \begin{tikzpicture}
        \matrix (m) [matrix of math nodes, row sep=3em, column sep=2.5em, text height=1.5ex, text depth=0.25ex]
        { \varinjlim \cU_c & U_c \\
        \varinjlim \cU_s  & U_s \\
        };
        \path[->,shorten >=0.5ex] (m-2-1) edge node[auto,font=\scriptsize] {$\varphi$} (m-1-1);
        \path[->>] (m-2-1) edge (m-2-2);
        \path[->] (m-1-1) edge node[auto,font=\scriptsize,inner sep=1pt] {$\cong$} (m-1-2);
        \path[->] (m-2-1) edge (m-1-2);
        \path[right hook->] (m-2-2) edge (m-1-2);
    \end{tikzpicture}
    \end{center}
    where the bottom row is surjective since $U_s=\langle U: U\in\cU_s\rangle$. It remains to show that $\varphi$ is injective. Then the diagonal map is an isomorphism onto the image of $U_s$ in $U_c$.

    It is well known that the direct limit of a system of groups can be constructed explicitly as the free product of the involved groups modulo the relations coming from the inclusions in the system of groups, see \cite[\S 1]{Ser:Tre:80}. Assume that $\varphi$ is not injective, and fix $1\neq g\in\varinjlim \cU_s$ with $\varphi(g)=1$.

    By the above description, the element $g$ can be written as a word $g=g_1\cdots g_n$ with $g_i\in U_i\in\cU_s$ and with $n>0$ minimal.

    If $n=1$, then $g\in U_1\in \cU_s\subset \cU_c$. But then $\varphi(g)=1$ implies $g=1$, since $\varphi$ is the identity on $U_1$.

    If $n>1$, then $\varphi(g) = \varphi(g_1)\cdots \varphi(g_n)=1$. By the description of the direct limit $\varinjlim \cU_c$, there is $1\leq i <n$ such that $\varphi(g_i)$ and $\varphi(g_{i+1})$ are contained in a common group $V \in \cU_c$, since the word $\varphi(g_1)\cdots\varphi(g_n)$ can be reduced to $1$.

    The group $V$ is either a root group or a root group interval. In both cases, there is a prenilpotent, non-nested pair of roots $(\alpha,\beta)$ with $c\in\alpha\cap\beta$ such that $\varphi(g_i)$ and $\varphi(g_{i+1})$ are contained in $U_{[\alpha,\beta]}$.

    But observe that $\varphi(g_i)$ and $\varphi(g_{i+1})$ are contained in $U_s$, so by Lemma \ref{lem:stabiliser_intersect_interval}, there is a group $U\in \cU_s$ such that
    \[
    \{\varphi(g_i), \varphi(g_{i+1})\} \subset U_{[\alpha,\beta]}\cap U_s \subseteq U
    \]
    Hence the product $g_ig_{i+1}$ can already be formed in the group $\varinjlim\cU_s$ and
    \[
    g_1\cdots g_{i-1}(g_ig_{i+1})g_{i+2}\cdots g_n
    \]
    is a representation of $g$ as a word of shorter length in $\varinjlim \cU_s$, which contradicts the minimality of $n$. So $\varphi$ is injective.
\end{Proof}

\noindent All of these characterisations allow us to give a different description of the Steinberg group.

\begin{Proposition}\label{p:tilde_g_steinberg}
    Assume that $G$ is 2-spherical and satisfies condition \eqref{condition}. The direct limit $\tilde G$ of the system of groups $U_s$ for all co-spherical simplices $s\in\Sigma$ with respect to their natural inclusions is isomorphic to the Steinberg group $\hat G$.
\end{Proposition}

\begin{Proof}
	The proof will be divided into two steps. In the first step, we construct a map $\hat G\rightarrow \tilde G$, the inverse map is constructed in the second step.
	\subparagraph{Step 1:} Of course, $U_\alpha\leq U_s$ for any co-spherical simplex $s$ such that $R_s\subset\alpha$. So we have an inclusion $U_\alpha \hookrightarrow U_s$ for any such simplex $s$. Roots are connected and chambers and panels are always co-spherical. So for any two co-spherical simplices $s,s'$ with $R_s,R_{s'} \subset \alpha$, there is a sequence of co-spherical simplices $s_1,\ldots,s_n$ such that
	\[
	U_s \alter{\subseteq}{\supseteq} U_{s_1} \alter{\subseteq}{\supseteq} \cdots \alter{\subseteq}{\supseteq} U_{s_n} \alter{\subseteq}{\supseteq} U_{s'},
	\]
	and the group $U_\alpha$ is contained in all of these. Hence the inclusion $U_\alpha\hookrightarrow \tilde G$, induced by any inclusion $U_\alpha\hookrightarrow U_s$ for any co-spherical simplex $s$ with $R_s\subset \alpha$, does not depend on the choice of $s$. A similar argument works for the inclusions $U_{[\alpha,\beta]}\hookrightarrow U_s$ for any co-spherical simplex $s$ such that $R_s\subset (\alpha \cap \beta)$.
	By the universal property of the colimit, there is hence a canonical homomorphism $\hat G \rightarrow \tilde G$ induced by these inclusions.

	\subparagraph{Step 2:} By Lemma \ref{lem:g_sinjection}, every group $U_s$ can be written as a direct limit of a subsystem of the direct system for $\hat G$, hence there is a canonical homomorphism $U_s\rightarrow \hat G$. Since these homomorphisms are all induced by the inclusions of root groups, they are compatible with the direct system for $\tilde G$. So there is a canonical epimorphism $\tilde G \rightarrow \hat G$.

	A closer inspection yields that the composition of these maps is by construction trivial on the root groups $U_\alpha$, hence they invert each other.
\end{Proof}

\noindent The importance of the Steinberg group for this paper comes from the following theorem by Caprace, for which we require some additional notation. For every root $\alpha$ write $X_\alpha\coloneq \langle U_\alpha\cap U_{-\alpha}\rangle$ and $H_\alpha=N_{X_\alpha}(U_\alpha)\cap N_{X_\alpha}(U_{-\alpha})$. For every pair of roots $(\alpha,\beta)$ let $X_{\alpha,\beta}\coloneq \langle X_\alpha\cap X_\beta\rangle$ and let $\Phi_{\alpha,\beta}$ be the rank 2 root system generated by $\alpha$ and $\beta$. Set $H_{\alpha,\beta}=\bigcap_{\gamma\in\Phi_{\alpha,\beta}} N_{X_{\alpha,\beta}}(U_\gamma)$. Then $(X_{\alpha,\beta},(U_\gamma)_{\gamma\in\Phi_{\alpha,\beta}},H_{\alpha,\beta})$ is a root datum of rank 2.

\begin{Theorem}[Caprace, Theorem 3.11 in \cite{Cap:O2S:07}]\label{th:caprace}
    Assume that the group $G$ is of $2$-spherical type and that its diagram has no isolated nodes. Suppose also that
	\begin{enumerate}
		\item For each $\alpha\in\Pi$, we have $[H_\alpha,U_\alpha]=U_\alpha$.
		\item For each pair of roots $(\alpha,\beta)$ in $\Pi$ such that the corresponding Coxeter group elements $s_\alpha$ and $s_\beta$ do not commute, the Steinberg group $\hat X_{\alpha,\beta}$ is the universal central extension of $X_{\alpha,\beta}$.
	\end{enumerate}
	Then $\hat G$ is the universal central extension of the little projective group $G^\dagger$.
\end{Theorem}

\noindent Condition \eqref{condition} is implied by the first condition, see the proof of Theorem 3.11 in \cite{Cap:O2S:07}. Furthermore, note that the second condition is essentially a condition on Moufang polygons. There is a complete classification of Moufang polygons satisfying this condition, see \cite{dMT:CE2:07}, stating that it is almost always valid except in the case of small fields. Using \cite[33.10-33.17]{TW:MP:02}, one can calculate that the first condition is also always fulfilled except for small fields.

\section{Homotopy groups of Wagoner complexes}\label{sec:homotopy}

There are various results to calculate the fundamental groups of complexes formed out of group cosets, for instance by Soulé, Brown and Abels-Holz. An equivalent point of view is to consider a group acting on a simplicial complex and to construct a presentation using stabilisers and the fundamental group of the quotient complex.

We require a very simple form of a result of this type. Unfortunately, we could not find such a simple formulation in the literature, so we include a precise statement of this result here. We then use the modern approach of complexes of groups as in \cite{BH:NPC:99} to give a simple proof.

\begin{Theorem}[II.12.20 in \cite{BH:NPC:99}]\label{th:homotopy_type_via_action}
    Let $X$ be a finite-dimensional simplicial complex and let a group $G$ act on $X$ simplicially. Let $Y\leq X$ be a subcomplex which is a strict fundamental domain for the action, that is, the complex $Y$ intersects each $G$-orbit in exactly one point.

    We assume that $Y$ is connected and simply connected, and we write $G_0 = \langle G_s : s\in Y\rangle$ for the subgroup generated by the stabilisers of simplices of $Y$. Then
    \begin{align*}
        \pi_0(\lvert X\rvert) &\cong G / G_0, \\
        \pi_1(\lvert X\rvert) &\cong \ker ( \varinjlim\{G_s: s\in Y\} \twoheadrightarrow G_0),
    \end{align*}
    where the direct limit is taken over the system of stabilisers $(G_s)_{s\in Y}$ with the natural inclusions.
\end{Theorem}

\noindent We will use the language of simple complexes of groups and of stratified spaces in the proof. Since the definition of these objects would enlarge this paper unnecessarily, we will try to treat these concepts in a `black box' fashion. All definitions and all the results we use can be found in \cite[Chapter II.12]{BH:NPC:99}.

\begin{Proof}
    All references given in this proof and all missing definitions are contained in Chapter II.12 in \cite{BH:NPC:99}. The main result we need for this proof is Proposition 12.20. We have adopted the notation of this proposition here. In the following, we will construct the objects appearing there and we will relate our situation to the situation in the proposition.

    The geometric realisation $\lvert X\rvert$ is a \emph{stratified simplicial complex} in the sense of Definition 12.1 and the discussion thereafter, where the strata are given by the $G$-images of closed simplices of $Y$. The set of strata $\cQ$ for the subcomplex $\lvert Y\rvert$ is given by the closed simplices in $Y$.
    
    The group $G$ then acts by strata preserving morphisms on $X$, and we construct the associated \emph{complex of groups} $G(\cQ)$ as in 12.15. By Proposition 12.20(2), we then have $\pi_0(\lvert X\rvert)\cong G/G_0$. If we consider the $G_0$-action on the connected component $G_0\cdot Y\subseteq X$, we obtain by Proposition 12.20(4) that 
    \[
    \pi_1(\lvert X\rvert) \cong \ker ( \widehat{G(\cQ)} \twoheadrightarrow G_0),
    \]
    where $\widehat{G(\cQ)}\cong\varinjlim\{G_s:s\in Y\}$ by Definition 12.12.
\end{Proof}

\noindent This can be used to calculate the number of connected components and the fundamental group of the geometric realisation of the Wagoner complex $|\cW(G)|$.

\begin{Theorem}\label{th:main_result} Assume that $G$ is a group of Kac-Moody type satisfying condition \eqref{condition} and such that the associated root system is $2$-spherical. Then we have
	\begin{align*}
        \pi_0(|\cW(G)|) &\cong G/G^\dagger. \\
        \intertext{If $G$ is of non-spherical type or of rank at least 3, we have}
		\pi_1(|\cW(G)|) &\cong \ker(\hat G\twoheadrightarrow G^\dagger).
    \end{align*}
\end{Theorem}

\begin{Proof}
    Consider the action of $G$ on the Wagoner complex $\cW(G)$. The stabiliser of a vertex $(U_s)$ is obviously the group $U_s$ itself, the stabiliser of a flag
    \[
    (U_{s_1}\subset U_{s_2} \subset \cdots \subset U_{s_r})
    \]
    is given by the smallest group $U_{s_1}$. The group $G_0 = \langle U_s : s\in\Sigma\rangle$ is equal to the little projective group $G^\dagger$. The direct limit of all stabilisers of simplices in $\cA$ with the natural inclusions is hence the group $\tilde G$.
    
    A strict fundamental domain for the action on $\lvert\cW(G)\rvert$ is given by the geometric realisation of the standard apartment $\lvert\cA\rvert$. If $G$ is of spherical type of rank $n$, then $\lvert\cA\rvert$ is homotopy equivalent to a $(n-1)$-sphere, and hence connected and simply connected if $n\geq 3$.

    If $G$ is of non-spherical type, then $\lvert\cA\rvert$ is homeomorphic to the Davis realisation of the Coxeter system associated to $G$ as defined in \cite[Chapter 7]{Dav:GTC:08}, which is contractible by Theorem 8.2.13 in \cite{Dav:GTC:08}.
    
    So we can apply Theorem \ref{th:homotopy_type_via_action} and use that we have $\tilde G \cong \hat G$ by Proposition \ref{p:tilde_g_steinberg}.
\end{Proof}

\begin{Corollary}\label{cor:main_result}
    In particular, if the prerequisites of Theorems \ref{th:caprace} and \ref{th:main_result} are satisfied, we have
        \[
	        \pi_1(|\cW(G)|) \cong H_2(G^\dagger;\Z).
        \]
\end{Corollary}

\section{The connection to algebraic \texorpdfstring{$K$}{K}-theory}\label{sec:ktheory}

Wagoner complexes were originally introduced by Wagoner in \cite{Wag:BSK:73} to provide a definition of higher algebraic $K$-theory. He then proved the following theorem, which is a special case of our result only if $R$ is a division ring.

\begin{Theorem}[Wagoner, Proposition 2 in \cite{Wag:BSK:73}]
	Let $R$ be a ring. If $n\geq 2$ then
	\[
    \pi_0(|\cW(\Gl_{n+1}(R))|) \cong \Gl_{n+1}(R)/\E_{n+1}(R) \cong H_1(\Gl_{n+1}(R);\Z) \cong K_1(R).
	\]
	If $n\geq 3$, then
	\[
    \pi_1(|\cW(\Gl_{n+1}(R))|) \cong \ker(\St_{n+1}(R)\twoheadrightarrow \E_{n+1}(R)) \cong H_2(\E_{n+1}(R);\Z) \cong K_2(R).
	\]
\end{Theorem}
In subsequent papers \cite{AKW:HAK:73} and \cite{AKW:HAK:77}, Anderson, Karoubi and Wagoner and in a final paper \cite{Wag:EAK:77} Wagoner proved that Wagoner's definition of higher algebraic $K$-theory coincides with Quillen's definition of $K$-theory in the sense that
\[
K_i(R) \cong \pi_{i-1}(|\cW(\Gl_{n+1}(R))|)\qquad\text{for }n\gg i, i>0.
\]

\noindent This leads naturally to the question whether
\[
\pi_{i-1}(|\cW(G)|) \cong \pi_i(BG^+)
\]
for any group with $\rank(G)\gg i$. This question would be interesting for hermitian $K$-theory, for example, by choosing the group $G$ to be a symplectic or orthogonal group. Unfortunately, we do not know whether these groups are isomorphic.

\section{Affine Wagoner complexes}\label{sec:affine}

For completeness, we also include the following construction, imitating the corresponding construction in \cite{Wag:HTp:75}. Let $G$ be a group with a spherical root datum with discrete valuation, compare \cite[Chapter 3]{Wei:SAB:09}. Let $\Delta$ be the associated affine building and let $\Sigma$ be a fixed apartment of $\Delta$. We will denote the set of spherical roots by $\Phi$ and the half-apartments of $\Sigma$ by $H_{\alpha,k}$ for $\alpha\in\Phi$, $k\in\Z$ with associated walls $\partial H_{\alpha,k}$.

\begin{Definition}
	Let $s$ be a simplex in $\Sigma$ and let $n \geq 1$. The \emph{$n$-cell around $s$} is defined to be the set
	\[
    C_n(s) \coloneq \bigcap_{\substack{\alpha\in\Phi,\, k\in\Z \\ R_s\subset H_{\alpha,nk}}} H_{\alpha,nk}.
	\]
\end{Definition}

\begin{Definition}
	For any $n\geq 1$ and any simplex $s\in\Sigma$, we write
	\[
	U^n_s \coloneq \langle U_{\alpha,k} : C_n(s)\subset H_{\alpha,k}, \alpha\in\Phi, k\in\Z \rangle.
	\]
	Analogously to Definition \ref{def:wagonercomplex}, for any $n\geq 1$, we define the \emph{affine Wagoner complex $\cW^n(G)$} to be the flag complex over the partially ordered set of cosets $\{ gU^n_s : s\in \Sigma, g\in G\}$.
\end{Definition}

\noindent For every fixed $n$, we have by construction $s \subseteq s'$ implying that $U^n_s \leq U^n_{s'}$. There are also simplicial epimorphisms from $\cW^n(G)$ onto the barycentric subdivisions of certain `coarsenings' of the affine building $\Delta$. For buildings of type $\tilde A_n$, this can be found in \cite[\S 5]{Wag:HTp:75}.

\begin{Remark}
    As above, we define the standard apartment of the Wagoner complex by
    \[
        \cA^n = \Flag(\{ gU^n_s : g\in G, s\in \Sigma \}).
    \]
    Obviously, the apartment $\cA^n$ is isomorphic to the flag complex over the set of $n$-cells $\{C_n(s):s\in\Sigma\}$. It is not hard to see that this complex is in turn isomorphic to the complex associated to the affine reflection group generated by the reflections at the walls $\partial H_{\alpha,nk}$ for $\alpha\in\Phi$ and $k\in\Z$. In particular, by rescaling by $\tfrac{1}{n}$, this complex is isomorphic to the Coxeter complex $\Sigma$, so $\cA^n$ is contractible.
\end{Remark}

\noindent In analogy to the construction above, we denote the direct limit of the system of groups $\{U^n_s : s\in\Sigma\}$ with their canonical inclusions by $\tilde G^n$. By the same argumentation as in the proof of Theorem \ref{th:main_result}, using that $\cA^n$ is contractible, we obtain

\begin{Proposition} For any $n\geq 1$, we have
	\[
		\pi_1(|\cW^n(G)|) \cong \ker( \tilde G^n \twoheadrightarrow G^\dagger)
	\]
\end{Proposition}

\paragraph{Construction} If $n$ divides $m$, then $C_n(s)\subseteq C_m(s)$ and hence $U^n_s \supseteq U^m_s$ for any $s\in\Sigma$. Consequently, the map $p^m_n: \cW^m(G) \rightarrow \cW^n(G)$ induced by $gU^m_s\mapsto gU^n_s$ is well-defined.

We obtain an inverse system of the simplicial complexes $\cW^n(G)$ with the maps $p^m_n$ for $n$ dividing $m$. This induces an inverse system on the fundamental groups $\pi_1(|\cW^n(G)|)$.

\begin{Definition}
	The inverse limit of this system is denoted by
	\[
	\pi_1^{\text{aff}}(G) \coloneq \varprojlim_n \pi_1(|\cW^n(G)|)
	\]
	and is called the \emph{affine fundamental group of $G$}.
\end{Definition}

\noindent The following theorem was already proved by Wagoner in \cite[\S 2]{Wag:HTp:75} for $G=\Sl_k(\F)$ over a local field $\F$ in a very different way.

\begin{Theorem}
	Let $\cW(G)$ be the Wagoner complex associated to the spherical building at infinity. Assume that the rank of the spherical building is at least three. Then, for every $n\geq 1$, there is a homomorphism
	\[
		\pi_1(|\cW(G)|)\rightarrow \pi_1(|\cW^n(G)|)
	\]
	that is compatible with the inverse system above. In particular, we obtain a homomorphism
	\[
	\pi_1(|\cW(G)|)\rightarrow \pi_1^{\text{aff}}(G).
	\]
\end{Theorem}

\begin{Proof}
	We will first construct a homomorphism $\hat G\rightarrow \tilde G^n$. By the universal property of the direct limit, it is enough to construct homomorphisms $U_\alpha$, $U_{[\alpha,\beta]}\rightarrow \tilde G^n$ which are compatible with the natural inclusions.

    Note that by the definition of root data with valuations, we have $U_\alpha= \bigcup_{k\in\Z}U_{\alpha,k}$, this is a filtration and hence a direct limit. Again, it is enough to construct a homomorphism $U_{\alpha,k}\rightarrow \tilde G^n$ for all $k\in \Z$. The group $U_{\alpha,k}$ is contained in all groups $U^n_s$ such that $C_n(s)\subseteq H_{\alpha,k}$. There are obviously such simplices, and we consider the inclusions $U_{\alpha,k}\rightarrow U^n_s$. By an argument similar to the one in the proof of Proposition \ref{p:tilde_g_steinberg}, using that half-apartments are connected and that all homomorphisms in the system of groups are inclusions, these maps are compatible with the system of groups for $\tilde G^n$, and we obtain homomorphisms $U_{\alpha,k}\hookrightarrow \tilde G^n$.

	For the groups $U_{[\alpha,\beta]}$, consider the filtration $U_{[\alpha,\beta],k}=\langle U_{\gamma,k} : \gamma\in[\alpha,\beta]\rangle$. Then the group $U_{[\alpha,\beta],k}$ is contained in all groups $U^n_s$ with $C_n(s)\subseteq \bigcap_{\gamma\in[\alpha,\beta]} H_{\gamma,k}$. This intersection contains a sector and hence such simplices exist. It is also connected and by the same argument as before, the inclusion is compatible with the system of groups and we obtain a homomorphism $U_{[\alpha,\beta]}\rightarrow \tilde G^n$.

    Since all of these maps are induced by inclusions, they are automatically compatible with the system of groups for $\hat G$. Hence we obtain a homomorphism $\hat G\rightarrow \tilde G^n$ which makes the following diagram commute:
	\[
\begin{xy}\xymatrix{
	1 \ar[r] & \pi_1(|\cW(G)|) \ar[r]\ar@.[d] & \hat G \ar[r]\ar[d] & G^\dagger \ar[r]\ar@{=}[d] & 1 \\
	1 \ar[r] & \pi_1(|\cW^n(G)|) \ar[r] & \tilde G^n \ar[r] & G^\dagger \ar[r] & 1.}
\end{xy}
	\]
	For standard reasons, the dotted map completing the above diagram exists. It is the desired homomorphism.
\end{Proof}

\noindent In the case of $G=\Sl_k(\F)$ over a local field $\F$, Wagoner proves in \cite{Wag:HTp:75}, that $\pi_1^{\text{aff}}(G)$ is isomorphic to the kernel of the universal central topological extension of $G=G^\dagger=\Sl_k(\F)$, sometimes denoted by $H_2^{\text{top}}(\Sl_k(\F))$. We do not know whether this is true in the general case. Wagoner's proof relies heavily on the fact that $H_2(\Sl_k(\F)) \cong H_2^{\text{top}}(\Sl_k(\F)) \oplus D$, where $D$ is infinitely divisible, and that $H_2^{\text{top}}(\Sl_k(\F))$ is isomorphic to the group of roots of unity of $\F$.

\nocite{Bro:CoG:82}
\bibliographystyle{alpha}
\bibliography{../../biblio}
\end{document}